\documentclass[12pt]{article}

\usepackage{amsfonts}
\usepackage{amssymb}
\usepackage{graphics}
\usepackage{amsmath, amsthm, latexsym}

\def \r{\mathbb{R}}
\def \p{\mathbb{P}}

\def \.{\cdot}

\def \dim{\tex\documentclass[12pt]{article}}

\usepackage{amsfonts}
\usepackage{amssymb}
\usepackage{graphics}
\usepackage{amsmath, amsthm, latexsym}

\def \r{\mathbb{R}}
\def \p{\mathbb{P}}

\theoremstyle{definition}

\newtheorem{Def}{Definition}[section]
\newtheorem{Rem}{Remark}[section]
\newtheorem{Prop}{Proposition}[section]
\newtheorem{Lem}{Lemma}[section]
\newtheorem{Th}{Theorem}[section]

\begin{document}

\title{On Rectification of Circles and an Extension of Beltrami's 
Theorem}
\author{Farzali Izadi}
\maketitle

{\it ``The only perfect geometrical figures are the straight 
line and the circle.''} 

\hspace{11.5cm} Plato

{\footnotesize {\bf Abstract.} 
The goal of this paper is to describe all local diffeomorphisms mapping
a family of circles, in an open subset of $\r^3$, into straight lines. This
paper contains two main results. The first is a complete description of the
rectifiable collection of circles in $\r^3$ passing through one point. 
It turns out that to be rectifiable all circles need to pass through some 
other common point. The second main result is a complete description of 
geometries in $\r^3$ in which all the geodesics are circles. This is a 
consequence of an extension of Beltrami's theorem by replacing 
straight lines with circles.}  

\vspace{0.5cm}

\vspace{0.5cm}

\vspace{1cm}

{\bf Introduction.}
The problem that the author solved has its origin in 
{\it Nomography}: \footnote{A discipline which was discovered by the 
French Civil Engineer, Docagne (1880). It turned out to have practical 
applications in many branches of science and technology 
(cf. \cite{hank,hilbert})} 
how to reduce a nomogram of aligned points to a circular 
nomogram? In more mathematical terms: 
what are local diffeomorphisms 
that send germs of lines to germs of circles? This question 
was initially posed for $2$-dimensional nomograms 
by G.S. Khovanskii 
and, solved by A.G. Khovanskii in that case (cf. ~\cite{kho}). 
Our result leads to a solution of the corresponding $3$-dimensional 
nomography. On the other hand, it is a continuation of M\"obius' 
classical work that describes all transformations taking lines to lines 
and circles to circles. It is also related to Beltrami's 
investigations. By Beltrami's classical theorem, 
all the geometries whose 
geodesics are locally straight lines have constant curvature (cf. 
~\cite{shar,wolf,carmo}). 
We prove that all geometries in $\r^3$ whose geodesics are locally 
circles must also have constant curvature. The similar fact 
in $\r^4$ is wrong. This  
was communicated to the author 
by V. Timorin (cf. \cite{timor} for details). 
We also give a complete description for all the metrics of these geometries.

The paper consists of 3 sections. In Section $1$, we give the precise 
definition of the rectification for a collection of circles passing through 
one point and we prove the fundamental theorem of rectification. 
This is our first main result. Section $2$ is devoted to the 
classification theorems for the \textit{rich families of circles}
i.e., families that look like the families 
of geodesics (they are 
point-wise rectifiable and for each point on each circle there 
is an open cone filled by the tangent lines to other circles). 
This classification gives rise to the following theorem: up to a 
projective transformation of the space of the image and a M\"obius 
transformation of the space of the inverse image, 
there exist exactly 
three diffeomorphism classes of rectifiable rich family of circles. 
In Section $3$, we discuss some applications of our results in 
Riemannian geometry. First we show how the rectification problem gives 
rise to an extension of Beltrami's theorem. More precisely, we will 
prove that if $U$ is a region in $\r^3$ such that all 
geodesics with respect to some metric ${\bf g}$ are circles, then 
${\bf g}$ has constant curvature. Then by using this fact, we will 
calculate all Riemannian metrics in which the geodesics are circles.     

\noindent{\bf Acknowledgments: }
I am grateful to my advisor Prof. Askold G. Khovanskii both for suggesting 
the problem and his help during the work. My thanks also go 
to Prof. Rick Sharpe for reading the manuscript and giving valuable 
suggestions, to Vladlen Timorin for providing me with the counterexamples 
for $4$-dimensional space and reading the first draft, and to Kiumars Kaveh for reading the final draft of the paper. 
Finally I gratefully acknowledge the support I received from the Iranian 
Ministry of Science and Technology and also from the Mathematics Department 
at the University of Toronto.

\section{Rectification of a Bundle of Circles in $\r^3$}
In this section, we study the behavior of the curves in a 
rectifiable bundle near the center of the bundle. Our main 
result in this section is a theorem which is called {\it ``the $54$-circles 
theorem.''} In order to prove this theorem, we first show that the 
coefficients of the Taylor polynomials of the curves in a central 
bundle are polynomial functions in the direction components of the 
tangent lines. We also show that these polynomials satisfy some 
symmetry relations. Assuming that all curves in the rectifiable 
bundle are circles, we prove that all polynomials are divisible 
by a specific irreducible polynomial. Using these divisibility conditions 
together with the above symmetric relations, the desired result can be 
easily obtained. 

Let us start with the following  
definition and notations.

\begin{Def}[Rectifiable Bundle of Curves]
A family of curves in the $3$-dimensional 
space $\r^3$, is called {\it rectifiable} 
near a point $p$ if there exists a neighborhood $U$ of $p$ and a diffeomorphism 
of $U$ taking all the curves in the family (more precisely, 
the portions of the curves contained in the region $U$) 
into straight lines (more precisely, into portions of lines lying in 
the image of the region $U$). Any collection of curves passing through the point 
$p$ is called {\it a bundle of curves} with center $p$. 
A bundle is called {\it simple} if distinct curves of the bundle 
have distinct tangent lines.
\end{Def}

If a bundle of curves with center at the point $p$ is rectifiable at $p$,
then it is simple. This means that each direction $(1,k,m)$, where 
$k,m \in \r^*=[-\infty,+\infty]$ corresponds to a unique curve 
$ \alpha_{(k,m)}$, where
$ \alpha_{(k,m)}$ is the intersection of the two surfaces 
$$F(x,y,z)=0 \qquad G(x,y,z)=0$$
and $(1,k,m)$ is the direction vector of the 
tangent line to the curve at the point $p$. 
For simplicity, we take the point $p$ to be at the origin. 
Now we are ready to state the following proposition which is useful 
in the proof of our main result.                 

\begin{Prop}\label{a}
Suppose that a simple bundle of curves  $\alpha_{(k,m)}$ 
is locally rectifiable near the origin by means of a class $ C^n$ 
diffeomorphism. Then for every $i, 1<i\le n$, there exist 
$2$-variable polynomials $ P_i$ and $ Q_i$ of degree at most $2i-1$ such that
$$y^{(i)}(0)=P_{2i-1}(k,m) \qquad z^{(i)}(0)=Q_{2i-1}(k,m)$$
where $y$ and $z$ 
are expressed in terms of $x$ on the curve  $\alpha_{(k,m)}$
\end{Prop}

To prove this proposition, we need a lemma 
that can be easily verified. It is in fact a sharpening of the 
{\it implicit function theorem} (cf. \cite{spivak}). 

\begin{Lem}\label{b}
Consider a curve $\alpha$ given by the equations $F(x,y,z)=G(x,y,z)=0$, 
where $F$ and $G$ are smooth functions on $\r^3$ such that
\begin{equation}
F(0,0,0)=G(0,0,0)=0,\hspace{.5cm} and \hspace{.5cm} d=det[\frac {\partial (F,G)}{\partial (y,z)}]_{(0,0,0)} \not =0. 
\end{equation}

Then by taking $x$ as a local parameter on $\alpha$ near $0$, 
the Taylor coefficients of $y$ and $z$, with respect to  
$x^i$'s are polynomials of degree $2i-1$ in the Taylor coefficients 
of $F$ and $G$.
\end{Lem}

\begin{proof}
These can be computed directly by using the implicit function theorem.
\end{proof}

\begin{proof}[Proof of the proposition] 
Consider a diffeomorphism $\Phi$ rectifying a bundle of curves 
$ \alpha_{(k,m)}(x)$. Let $T$ be an arbitrary nonsingular linear 
mapping of the space. The diffeomorphism $ T \circ \Phi$ also 
rectifies the bundle $\alpha_{(k,m)}(x)$. By a suitable choice 
of $T$ we can arrange
that the rectifying diffeomorphism has the identity differential 
at the origin. The rectifying diffeomorphism is now given by the map:
$$T \circ \Phi=(f,g,h)=(u,v,w)$$
where

\begin{equation}
f=x+ \cdots,\\ 
g=y+ \cdots,\\
h=z+ \cdots.
\end{equation}

and $\cdots$ denotes the higher terms in the expansion of each function. 
In the $uvw$-space the bundle of curves $ \alpha_{(k,m)}(x)$ 
is given by the equations $ v=k\.u,w=m\.u$. Consequently, the curves $ 
\alpha_{(k,m)}(x)$ are given by the equations: 
$F(x,y,z)=0, G(x,y,z)=0$, where, $F=g-kf, G=h-mf.$ Now the 
coefficients $ a_{pqr}$ and $ b_{pqr}$ in the Taylor polynomial 
of the functions  $ F$  and  $ G$  depend linearly on 
$\{ k,a_{010}=1\}$ and on $\{ m,b_{001}=1\}$ respectively. 
The assertion now follows from Lemma \ref{b}. 
\end{proof}

The next proposition is the second useful fact concerning our main result. 
\begin{Prop}\label{c}
Under the same assumptions as of Proposition \ref{a}, there exist $3 \times 3$ 
symmetric matrices $A, B$, and $C$ such that 
$$y^{(2)}(0)=\langle(B+kA) \lambda,\lambda \rangle, \hspace{.5cm}
z^{(2)}(0)=\langle(C+mA) \lambda, \lambda \rangle$$
where $\lambda=(1,k,m)$ is the tangent vector at the origin.
\end{Prop}

\begin{proof}
In order to prove this, suppose that the curve 
$ \alpha_{(k,m)}(x)=(x,y(x),z(x))$ lies on the surface 
$ F(x,y,z)=0$, where $F$ is any smooth function. So 
$ F(\alpha_{(k,m)}(x))=0$. For simplicity we write 
$\alpha$ in stead of $\alpha_{(k,m)}$. 
By differentiating both sides of the equation we get
$$ \langle \nabla F(\alpha(x)),\alpha^{'}(x) \rangle=0.$$ 
which implies that
$${\langle \frac{d}{dx} \nabla F(\alpha(x)),\alpha^{'}(x)\rangle+
\langle \nabla F(\alpha(x)),\alpha^{(2)}(x)\rangle=0}. \qquad (*)$$
where
$${\frac {d}{dx} \nabla F(\alpha(x))=H(\alpha(x))\alpha^{'}(x)=
(\frac{\partial ^2F}{\partial x_i \partial x_j})\alpha^{'}(x)}.$$
This expression holds for any surface $F$ containing the graph of the 
curve $\alpha(x)$, in particular for $F=v-ku=g-kf$, where, $f$ and $g$ 
are the same functions in (2). Since 
$$\nabla F(\alpha(0))=(-k,1,0), \hspace{.5cm} 
\alpha^{(2)}(0)=(0,y^{(2)}(0),z^{(2)}(0))$$
and $H(\alpha(0))=B+kA$, where $B=
(\frac{\partial ^2g}{\partial x_i \partial x_j})(0)$ and 
$A=(\frac{\partial ^2(-f)}{\partial x_i \partial x_j})(0)$.

By (*) we get
$$ y^{(2)}(0)=\langle (-k,1,0),\alpha^{(2)}(0)\rangle=
-\lambda H \alpha((0)){\lambda}^T$$ 
where, $\lambda=(1,k,m)$. Hence
$$y^{(2)}(0)=-\langle(B+kA)\lambda,\lambda \rangle.$$

Similarly, for calculating $ z^{(2)}(0)$, we let $G=w-m\.u=h-m\.f$, 
where $f$ and $g$ 
are the same functions as in $(2)$. 
Clearly, $ \nabla G_m(\alpha(0))=(-m,0,1)$, and
$H(\alpha(0))=C+mA$, where $C=(\frac{\partial^2h}
{\partial x_i \partial x_j})(0)$. 
Hence $z^{(2)}(0)=-\langle(C+mA) \lambda,\lambda \rangle.$ 
\end{proof}

\begin{Rem} 
These propositions  can be easily extended to arbitrary dimensions, 
but we don't need to do this. 
\end{Rem}

Before stating the first main result,
we give a definition.
\begin{Def}
Let us say that $54$ lines passing through the origin 
are generic if there exists a unique homogeneous cone of degree $9$ 
containing them all. In this case, the corresponding $54$ points in 
$\p^2$ are called \textit{9-good}.
\end{Def}

One can easily show that for almost all $54$ lines passing 
through the origin, the above condition holds. Having said this, we 
have the following: 

\begin{Th}[The $54$-circles theorem]\label{d} 
Consider a simple
bundle of circles passing through the origin such that the set of 
tangent lines of the circles contains a $54$ generic lines, then 
there exists a local diffeomorphism about the origin mapping all 
the circles into straight lines if and only if all the circles in 
the bundle pass through one common point distinct from the origin.
\end{Th}

\begin{proof}
In one direction the proof is obvious. 
In fact, if the bundle passes through the second point $Q$, 
then we can make an inversion with respect to a sphere 
centered at this point. But in the opposite direction, 
it is rather complicated and follows as:

The bundle of circles passing through the origin 
in $\r^3$ can be written explicitly by a system of equations
consisting of two spheres. The simplicity condition easily 
implies that this system depends only on two parameters, 
namely the components of tangent vector at the origin. 
Using this fact, the above system of equations can be 
expressed in the following form

\begin{equation}
\left \{
\begin{array}{c}
 y=kx+A(x^2+y^2+z^2)\\
 z=mx+B(x^2+y^2+z^2)
\end{array}
\right.
\end{equation}
where, $ A=A(k,m)$ and $ B=B(k,m)$ are 
some functions of the parameters
$k$ and $m$. We show that the rectifiability of 
the bundle is equivalent to the linearity of the 
functions  $A$ and $B$. We wish to solve the 
equations for the circles in the bundle up to 
terms of fourth order. By 
differentiating both equations with respect to $x$ for 
the Taylor series of $y(x)$ and $z(x)$ we obtain

$$y(x)=kx+\phi_2(k,m)x^2+\phi_3(k,m)x^3+\phi_4(k,m)x^4+\cdots$$
$$z(x)=mx+\psi_2(k,m)x^2+\psi_3(k,m)x^3+\psi_4(k,m)x^4+\cdots$$
where, by letting $f=1+k^2+m^2$, we have

\begin{equation}
\begin{array}{c}
 \phi_2=Af\\
 \psi_2=Bf\\
\end{array}
\end{equation}

\begin{equation}
\begin{array}{c}
 \phi_3=2A(kA+mB)f\\
 \psi_3=2B(kA+mB)f\\
\end{array}
\end{equation}

\begin{equation}
\begin{array}{c}
\phi_4=A(A^2+B^2)f^2+4A(kA+mB)^2f\\
 
\psi_4=B(A^2+B^2)f^2+4B(kA+mB)^2f\\
\end{array}
\end{equation}
where $\phi_l$ and $\psi_l$ are polynomials of 
degrees at most $2l-1$ in the two variables $k$ and $m$ 
(Proposition \ref{a}), and $A, B$ are, at the outset, just 
rational functions in $k$ and $m$. We want to show that 
$A$ and $B$ are polynomials of degree 1.

Multiplying equations (5) by $f$ yields

\begin{equation}
\begin{array}{c}
 f\phi_3 = 2Af(kA+mB)f = \phi_2(k\phi_2+m\psi_2) \\
 f\psi_3 = 2Bf(kA+mB)f = \psi_2(k\phi_2+m\psi_2)\\
\end{array}
\end{equation}
According to Proposition \ref{a}, the functions 
$ (\phi_2,\psi_2) ,(\phi_3,\psi_3)$ and $ (\phi_4,\psi_4)$
are polynomials of degree at most 3,5, and 7 in $k,m$ 
respectively. Equations (7) are satisfied for all values
of $(1,k,m)$ corresponding to circles in the bundle, 
i.e., by at least a 54 directional vectors corresponding
to 54 generic lines. Since 2-variable polynomials of 
degree 5 which coincide at 54 directional vectors, 
coincide identically, so these equations are in fact
identities. These identities imply either: $f$ 
divides $(k\phi_2+m\psi_2)$, or: $f$ divides both 
$\phi_2$ and $\psi_2$. In the latter case equations 
(4) show that $A$ and $B$ are polynomials of degree 1. 
So we may assume that $fg=k\phi_2+m\psi_2$ for some polynomial $g$.

Now
\begin{eqnarray*}
fg &=& k\phi_2 + m\psi_2 \cr
&=& kAf + mBf \qquad \textup{(by equation (4))} \cr
&=& (kA + mB)f \cr
\end{eqnarray*}
Thus $kA + mB = g$ is a polynomial. 

Multiplying equations (6) by $f$ yields
\begin{eqnarray*}
f\phi_4 &=& (Af)((Af)^2 + (Bf)^2 + 4(Af)(kA + mB)^2f \cr
&=& \phi_2(\phi_2^2 + \psi_2^2) + 4\phi_2g^2f \cr
\end{eqnarray*}
Similarly    
$$f\psi_4 = \psi_2(\phi_2^2 + \psi_2^2) + 4\psi_2g^2f$$

By the same reasoning as before, these equations are again 
identities. This shows that either: $f$ divides $(\phi_2^2 + \psi_2^2)$,
or: $f$ divides both $\phi_2$ and $\psi_2$. In the second case, 
as before, we are done, so we may assume that $f|(\phi_2^2 + \psi_2^2)$.

Equation (5) gives
\begin{eqnarray*}
m\phi_3 + k\psi_3 &=& 2mA(kA = mB)f + 2kB(kA = mB)f \cr
&=& 2f(mA + kB)(kA + mB) \cr
&=& 2f(km(A^2 + B^2) = 2(k^2 + m^2)AB) \cr
\end{eqnarray*}

So $f(m\phi_3 + k\psi_3) = 2km(\phi_2^2 + \psi_2^2) + 2(k^2 + m^2)\phi_2\psi_2.$
Since we already know that $f|(\phi_2^2 + \psi_2^2)$ 
it follows that $f|\phi_2\psi_2$. Thus we also have 
$f|(\phi_2 \pm \psi_2)^2$ and thus $f|\phi_2 \pm \psi_2$ and finally 
$f|\phi_2$ and $f|\psi_2$. 
Since these polynomials are polynomials of degree 3 in $k$ and 
$m$, $A(k,m)$ and $B(k,m)$ would be polynomials of degree 1, in 
$k$ and $m$ ,i.e. 
$$ A(k,m)=ak+bm+c, \qquad B(k,m)=dk+em+f.$$ 
Now by using the symmetric relations in Proposition \ref{c}  
we see that $b=d=0$ and $a=e$. Hence the functions $A$ and $B$ 
are in the form
$$A(k,m)=\alpha k+\beta, \hspace{.5cm}B(k,m)=\alpha m+\gamma.$$
These linear functions show that our rectifiable bundle 
of circles necessarily has the form
\begin{equation}
\left \{
\begin{array}{c}
 S_1+kS_2=0\\
 S_1+mS_3=0
\end{array}\
\right.
\end{equation}
where, $ S_1=0,S_2=0$ and $ S_3=0$ are the 
equations for certain non-tangent spheres passing 
through the point  $(0,0,0)$. We denote by $Q$ the 
second point of intersection of the spheres $ S_1=0,S_2=0$ 
and $ S_3=0$. All the circles in the bundle pass through $Q$.
In order to rectify such a bundle of circles, it suffices 
to map the point $Q$ to infinity via an inversion. 
Now the proof of theorem is complete.  
\end{proof}

\begin{Rem}
The analogous result fails 
in $\r^4$. Suppose that this is not the case. Then for
every point $p$ the bundle of all circles (geodesics)
passing through the point $p$, being rectifiable by 
the exponential map, should pass through some other
common point. Now it can be easily shown 
(as in the proof of Theorem \ref{d}) that there exists a germ 
of local diffeomorphim mapping all the circles 
in a neighborhood into staight lines. This would 
give rise to a $4$-dimensional extension of Beltrami's
theorem which in turn implies that the corresponding 
metric in $\r^4$ has constant curvature. 
But there is a famous example of a metric in $\r^4$ i.e. 
{\it Fubini-Study} metric 
which has circle geodesics 
but non-constant curvature (cf. ~\cite{grif} or ~\cite{arnold} 
for details).
\end{Rem}

\section {Classification Theorems}
Consider the space S of equations of spheres in $\r^3$ i.e., 
the space of non-zero polynomials of the form
$$V=\{a(\sum_{i=1}^3x_i^2)+\langle b,x \rangle+c 
|a,c \in \r, b \in \r^3\}.$$

Clearly, every element of this form is 
defined up to a factor. Thus the space $V$ is  
isomorphic to the projective space $\r P^4$. A 
projective subspace $L$ of $V$ of dimension $k (k=1,2,3)$ 
is called a $k$-dimensional linear system of spheres. 
Among all different linear systems, there are three 
systems which are closely related to the three 
geometries of Lobachevski, Euclid, and Riemann: the 
linear system of all spheres orthogonal, respectively 
to a fixed sphere of positive radius: $\sum_{i=1}^{3}x_i^2=1$, 
zero radius: $\sum_{i=1}^{3}x_i^2=0$, and imaginary radius: 
$\sum_{i=1}^{3}x_i^2=-1$ . These three different linear 
systems can be expressed as:

\begin{enumerate}
\item $A(\sum_{i=1}^3x_i^2)+\langle B,x \rangle+A=0.$
\item $\langle B,x \rangle+D=0.$ 
\item $A(\sum_{i=1}^3x_i^2)+ \langle B,x \rangle-A=0.$
\end{enumerate}

\begin{Def}
A $3$-dimensional \textit{net of spheres} 
is any set of spheres, the equations of which lie in some 
$3$-dimensional linear system but not in any $2$-dimensional 
linear system.
\end{Def}

\begin{Def}[Characteristic map]
A \textit{characteristic map} of $3$-dimensional net is a map 
$\Phi:\r^3 \longrightarrow \r P^3$
defined by
$$\Phi(X)=[S_1(X): S_2(X): S_3(X) : S_4(X)]$$
where, $S_1,S_2,S_3,S_4$ are any 
$4$ independent quadratic polynomials in the space $V$. 
A characteristic map $\Phi$ depends on the choice of 
the polynomials $S_i$ and is therefore defined up to 
a projective transformation.
\end{Def}
\begin{Def}[Degenerate point] 
The point $(x_1,x_2,x_3)$ of a characteristic map 
$\Phi$ is called a \textit{degenerate point} of $3$-dimensional net of 
spheres if $\Phi$ has the zero Jacobian at that point. The 
degenerate points of the three linear systems of spheres 
indicated above consist, respectively, of the points on 
the sphere $\sum_{i=1}^{3}x_i^2=1$, the point $(0,0,0)$, and the 
empty set. 
\end{Def}
\begin{Def}[Rich family of circles in $\r^3$]. 
Let $\Delta$ be a family of circles in some domain $U$. 
$\Delta$ is called a \textit{rich family}, if there exists a subfamily 
$\Gamma \subseteq \Delta$ such that

\begin{enumerate}
\item For each $P \in U$ there exists a circle 
$\gamma \in \Gamma$ such that $P \in \gamma$.
\item If $\gamma \in \Gamma$ and $P \in \gamma$, 
then there exists an open cone $K_P$ -It is assumed that 
the cone depends continuously on the point P- such that 
the tangent line of $\gamma$ at the point $P$ lies inside $K_P$, 
and any other direction in $K_P$ corresponds to a circle in $\Gamma$.
\end{enumerate}
\end{Def}

Next two theorems give a complete description of all 
local diffeomorphisms which rectify a rich family of circles 
in a domain $U$. 

\begin{Th} \label{g}
A rich family of circles in $\r^3$ 
in a neighborhood of the point $P$ is rectifiable if and only
if there exists a germ of a diffeomorphism 
$\Phi:(\r^3,P) \rightarrow \r \p^3$ 
given by
$$\Phi(X)=[S_1(X):S_2(X):S_3(X):S_4(X)].$$
where 
$$S_i(X)=a_i(x^2+y^2+z^2)+b_ix+c_iy+d_iz+e_i \qquad (i=1,\ldots,4).$$
with a non-zero Jacobian such that every circle in the 
family is the inverse image of a line under $\Phi$.
\end{Th}

\begin{proof}
Let us first suppose that the
rich family of circles is rectifiable. Let $c=0$ be the equation
for some circle in the family passing through the point $P$, 
with a tangent lying inside the cone $ K_p$. Let $A$ and $B$ 
be two points on the circle $c=0$ lying close to $P$ but on 
different sides. The circles in the family passing through 
the points $A$ and $B$ form rectifiable bundles by our assumption.
Hence by theorem \ref{d} they pass through the points $C$ and 
$D$ respectively distinct from $A$ and $B$. Through each 
point $Q$ close to the $P$, (i.e., $Q$ contained in both 
$K_A$ and $K_B$), there exist circles in the family of 
circles passing through the points $A$, $C$ and in the 
family of circles passing through the points $B$, $D$. 
Again by theorem \ref{d} all circles in the family passing 
through the point $Q$ pass a single point $S$ distinct 
from the point $Q$. Now suppose that the lines containing 
the segments $\overline{AC}$ and $\overline{BD}$ intersect 
at some point $O$. According to the different positions of 
the point $O$, we have the following different cases. 

\noindent {\it Case (1)}. The point $O$ lies outside $c$. By 
definition of the power of the point $O$ with respect to the 
circle $c$, we have: $\overline{OA}\cdot\overline{OC}=
\overline{OB}\cdot\overline{OD}$. Let us denote this number 
by $r^2$. Let $S$ be a sphere of radius $r$ centered at the 
point $O$. It is easy to see that $R^2+r^2=\overline{OM}^2$, 
where $R$ and $M$ are the radius and center of the circle $c$ 
respectively.

\noindent {\it Case (2)}. The point $O$ is on the circle $c$. 
In this case, the sphere $S$ is just a single point $O$. In 
other words, we have a sphere of radius zero.

\noindent {\it Case (3)}. The point $O$ lies inside $c$. 
In this case, one can easily show that $R^2-r^2=\overline{OM}^2$
, where $r,R,O$ and $M$ are the same as the case(1).

Now by mapping the spheres in cases (1) and (3) to the unit
sphere at the origin,
and similarly the point $O$ in case (2) into origin, we can
easily see that the equations of circles lie on the following 
three different spheres:
$$Ax+By+Cz+D(1+x^2+y^2+z^2)=0,$$ 
$$Ax+By+Cz+D(0+x^2+y^2+z^2)=0,$$
$$Ax+By+Cz+D(1+x^2+y^2+z^2)=0.$$
To complete the proof we only need to show that for 
each rich family of rectifiable circles the characteristic map 
is the same from one point to another. To prove this last assertion,
suppose that $P, Q \in U$ such that $P \ne Q$. Since any curve 
joining the two points $P$ and $Q$ is a finite curve or a segment,
we can cover this curve by a finite number of balls such that the
center of each ball lies in the next. Each characteristic map is 
a analytic function and any two characteristic maps coincide in 
the intersection of their domain, so they coincide identically 
by the theorem of analytic continuation.

Thus the equations for all the circles in our rectifiable family 
lie in two different linear combinations of the equations of spheres
(including planes) $ S_1=0, S_2=0, S_3=0$ and $ S_4=0$. Moreover, 
we can easily see that near a nondegenerate point of a rich family
there does not exist any rectifiable rich subfamily (this can be 
verified separately for the three linear systems of circles). 
Near a nondegenerate point of the family, the family is rectified
by the characteristic transformation
$$\Phi(x,y,z)=[S_1(x,y,z):S_2(x,y,z):S_3(x,y,z):S_4(x,y,z)].$$
where 
$$S_i(X)=a_i(x^2+y^2+z^2)+b_ix+c_iy+d_iz+e_i, (i=1,\ldots,4).$$
\end{proof}

It remains to show that up to a projective transformations there 
exists no other rectifying map. This is an immediate consequence 
of the following Lemma.
  
\begin{Lem}\label{e}
A local diffeomorphism of the space which 
sends a rich family of lines into lines is a projective transformation.
\end{Lem}

\begin{proof}
It is well-known that a 
homeomorphism which sends all lines into lines is a projective 
transformation. The proof of this fact based on constructing an 
everywhere M\"obius flat net (cf. \cite{klein}).
To this end, suppose that four lines $l_i (i=1,\ldots,4)$
in a plane $\Pi$ are in general position. Suppose $F$ is a map 
which sends these four lines into another four lines $F(l_i) 
(i=1,\ldots,4)$ in general position. Now for any projective 
transformation $T$, $ T(l_i) (i=1,\ldots,4)$ are also in general
position. Then there exists a projective transformation $U$ 
such that $U$ maps $T(l_i)$ into $F(l_i)$. So without loss of 
generality, we may assume that $T(l_i)=F(l_i) (i=1,\ldots,4)$. 
For each $i (i=1,\ldots,4)$ let us denote this line by $m_i$. 
Suppose that $A, B, C, D, $ and $a, b, c, d, $ are four vertices 
of the quadrilateral formed by $l_i$ and $m_i$, $(i=1,\ldots,4)$ 
respectively. Since both maps $F$ and $T$ map the diagonals of
the first quadrilateral to the diagonals of the second one, 
it follows that $F(P)=T(P)$, where $P$ and $F(P)$ are the 
intersection points of the pairs of the diagonals respectively.
Continuing in this way, we get a countable dense subset 
of some neighborhood such that $F=T$ on this subset. Since both $F$
and $T$ are continuous, we have $F=T$ on this neighborhood. Now 
the proof of lemma \ref{e} based on the same argument. If $P$ is a 
point in the neighborhood $U$, then we apply the same reasoning
for any plane $\Pi$ passing through the point $P$. First of all,
we can map any four lines in general position into a parallelogram
via a projective transformation. Secondly, by definition of a rich
family of lines for any point $a \in U$, there exists a line $l$ 
passing through the point $a$. Since $a \in l$, there exists a
cone $K_a$ such that $l$ lies inside $K_a$. For any other 
point $b \in l$, there exists another cone $K_b$ such that 
$l$ lies inside $K_b$. Now we can easily construct a 
parallelogram such that all four sides as well as its 
diameters contained in our rich family. Finally, we 
construct a {\it M\"obius flat net} inside this parallelogram,
all of whose lines lie in the rich family. This implies that
the mapping $F$ is locally projective. The connectedness of
the region $U$ now implies that $F$ is projective.
\end{proof}

We proved that

\begin{Th}\label{f}
Up to a projective transformation of 
the space of the image and a M\"obius transformation of the space
of the inverse image, there exist exactly three local 
diffeomorphisms which rectify rich families of circles. 
They are given by:
\begin{enumerate}
\item  $\Phi(x,y,z)=[x:y:z:1-x^2-y^2-z^2]$
\item  $\Phi(x,y,z)=[x:y:z:1]$
\item  $\Phi(x,y,z)=[x:y:z:1+x^2+y^2+z^2]$   
\end{enumerate}
\end{Th}

\section{Applications in Riemannian Geometry}

As is well-known, there are three classical geometries 
in which the geodesics in some local coordinate system
are straight lines. Beltrami's theorem together with 
the Minding-Riemann theorem, both of which are n-dimensional
results, ensure that these geometries are the only ones with
these properties. Fixing $n=3$, we will show that there are
precisely three classical geometries whose geodesics are circles.
This will be a consequence of the $3$-dimensional extension
of Beltrami's theorem which we are going to prove by replacing
straight lines with circles. However, this result, unlike the
previous statement, does not hold in arbitrary dimension. 
In fact, there is a very natural metric in $\r^4$ whose 
geodesics are circles, but it does not have constant curvature and
hence it does not coincide with the three classical geometries.

\subsection{Riemannian Geometry}
In this subsection, 
we give a proof of the $3$-dimensional extension 
of the Beltrami theorem. Next, we recall the three 
metrics in which the geodesics are straight lines. 
In the end, we use these metrics to calculate 
the corresponding metrics with circle geodesics. 

\begin{Th}[ $3$-dimensional extension of the Beltrami theorem for circles]
Let $U$ be a region in $\r^3$ and $g$ a Riemannian metric on $U$.
If all geodesics in $U$ are circles (or parts of cirlce), then $g$ has
costant curvature.
\end{Th}

\begin{proof}
First of all, note that all geodesics (circles) passing through 
one point are rectifiable by the exponential map. 
Secondly, by Theorem \ref{d}  any rectifiable bundle in $\r^3$
passes through some other common point distinct from the
center. On the other hand, all geodesics in a region 
$U$ form a rich family of circles. Now by using Theorem \ref{g} 
together with Theorem \ref{f} we get the three different 
characteristic maps. By Theorem \ref{f}, every characteristic 
map is a geodesic map (mapping the geodesics of the first 
space to the geodesics of the second). The Beltrami 
theorem now implies that $U$ has constant curvature $k$. 
By the Minding-Riemann theorem we see that $U$ 
isometric to a part of the elliptic space if $k>0$, Euclidean 
if $k=0$, and hyperbolic if $k<0$.
\end{proof}

\begin{Rem}
First of all, this result is
absolutely different from that of Beltrami's. In the 
Beltrami theorem the geodesics are already given  and
we look for the corresponding rectification. But here
we  deal with a $6$-dimensional family of circles in 
${\bf \r^3}$ and we wish to choose among this family
all $4$-dimensional families such that they could 
be families of geodesics. The problem of choosing 
these families relies heavily on the rectification
problem. 

Secondly, the Beltrami theorem holds in
any dimension, but our result can not be 
extended to dimension $4$. As we have already mentioned 
some very natural examples of Riemannian metrics in 
dimension $4$ in which the geodesics are circles but 
the curvature is not constant, as we have already mentioned, 
is the {\it Fubini-Study} metric.

To have a view of the nature of the metric $g$ in 
the space of rectifiable families of circles we need to 
draw our attention to some models of geometries in which the 
geodesics are straight lines (cf. \cite{rat}).
\end{Rem}

\begin{Rem}
The following discussion is also 
true for $n=2$. Clearly for $k=0$, the model is the Euclidean 
space $\r^3$ with the corresponding Riemannian metric
$$ds^2=dz_1^2+dz_2^2+dz_3^2.$$
Hence by the affine characteristic map $\Phi(z)=x$ we have
$$ds^2=\sum_{i=1}^3dx_i^2.$$
For $k<0$, we use the gnomonic projection of $D^3$ onto $H^3$,
where 
$$D^3=\{z\in \r^3:\sum_{i=1}^3z_i^2 <1\},$$ 
and
$H^3=\{z\in \r^4:z_1^2+z_2^2+z_3^2-z_4^2=-1, z_4>0\}.$
Identify $\r^3$ with $\r^3\times \{0\}$ in $\r^4$. 
The gnomonic projection $\mu$ of $D^3$ onto $H^3$ is 
defined to be the composition of vertical translation 
of $D^3$ by $e_4$ followed by the radial projection to 
$H^3$. An explicit formula for $\mu$ is given by
$$\mu(y)=\frac{y+e_4}{\Vert y+e_4\Vert}$$
where 
$$\Vert y+e_4 \Vert ^2=1-\vert y \vert ^2=1-\sum_{i=1}^3 y_i^2.$$
First we note that the element of hyperbolic arc length of $H^3$ is
$$\Vert dz \Vert ^2=\sum_{i=1}^3 dz_i^2-dz_4^2.$$
If $y.dy$ denotes the standard inner product 
between $y=(y_1,y_2,y_3)$ and $dy=(dy_1,dy_2,dy_3)$, 
then
\begin{equation}
\vert dz \vert ^2=\sum_{i=1}^3 dz_i^2-dz_4^2=
\frac{(1-\vert y \vert ^2)\vert dy \vert^2+(y.dy)^2}{(1- \vert y \vert ^2)^2}.
\end{equation}
For $k>0$, we similarly use the gnomonic
projection of $\r^3$ onto $S^3$, the unit sphere. 
To do this, we identify $\r^3$ with $\r^3 \times \{0\}$ in  
$\r^4$. Then the gnomonic projection
$$\nu:\r^3\rightarrow S^3$$
is defined to be the composition of the vertical 
translation of $\r^3$ by $e_4$ followed by the radial 
projection to $S^3$. An explicit formula for $\nu$ is given by    
$$\nu(y)=\frac{y+e_4}{\vert y+e_4 \vert}$$
where $\vert y+e_4 \vert$ is the Euclidean 
norm of $y+e_4$.
Since the element of spherical arc length of $S^3$ is 
the element of Euclidean arc length of $\r^4$ restricted 
to $S^3$, the arc length $ds$ of $S^3$ is given by
$$ds^2=\sum_{i=1}^4dz_i^2.$$
Thus
\begin{equation}
ds^2=\sum_{i=1}^4dz_i^2=\frac{\sum_{i=1}^3dy_i^2}{1+\vert y \vert^2}
-\frac{(y.dy)^2}{(1+\vert y \vert^2)^2}.
\end{equation}
\noindent This gives rise to an elliptic model in which
all geodesics are straight lines in $\r^3$.

Now we are ready to calculate the Riemannian metrics 
of our three geometries in which the geodesics at 
every point are circles.

Clearly for $k=0$ this is the standard
metric
$$ds^2=\sum_{i=1}^3 dx_i^2.$$
For the hyperbolic case, we set $k=-1$, 
and use the affine characteristic map $\Phi:\r^3 \rightarrow \r^3$ 
defined by
$$\Phi(x)=\frac{x}{1+\vert x \vert ^2}=y$$
where
$$\vert x \vert^2=\sum_{i=1}^3 x_i^2.$$
Let $w=({1+\vert x \vert^2})^{-1}$. 
Then $y_i=x_iw$, and $dy_i=dx_iw-2x_i(x.dx)w^2 \hspace{.5cm} i=1,2,3.$

Substituting these expressions into 
hyperbolic metric with geodesics as straight 
lines namely in (4.1) we get the following 
metric with geodesics as circles.
$$ds^2=\frac{1}{1+\vert x \vert^2+\vert x \vert^4}
(\vert dx \vert^2-\frac{3(x.dx)^2}{1+\vert x \vert^2+ \vert x \vert^4}).$$

Finally for the elliptic case we set $k=1$, and we 
use the affine characteristic map $\Phi: \r^3 \rightarrow \r^3$ defined by
$$\Phi(x)=\frac{x}{1-\vert x \vert^2}=y.$$
Let $t=({1-\vert x \vert^2})^{-1}$. 
Then $y_i=x_it$, and $dy_i=dx_it+2x_i(x.dx)t^2, \hspace{.5cm}i=1,2,3.$
Substituting these expressions into 
elliptic metric with geodesics as straight 
lines namely in (4.2) we get the following metric 
with geodesics as circles.
$$ds^2=\frac{1}{1-\vert x \vert^2+\vert x \vert^4}
(\vert dx \vert^2+\frac{3(x.dx)^2}{1-\vert x \vert^2+ \vert x \vert^4}).$$
\end{Rem}

\subsection{Open questions}
In the end, I state some open problems posed by A.G. Khovanskii. 
This work opens a large field for further investigations. 
Some most obvious questions to be asked are as follows:
Firstly, what is going on in higher dimension (both of 
the above main results are wrong in $\r^4$ (for details 
see the recent paper by V. Timorin \cite{timor}), in fact there are 
at least five different geometries having circle geodesics 
in $\r^4$.) 
Secondly, are there some results in the same spirit for more 
general classes of curves (say, algebraic curves of given 
degree with fixed leading term - the existence of a unique 
asymptotic cone seems to be important)?

\end{document}